\documentclass[12pt,a4paper]{article}
\usepackage{amsmath}
\usepackage{apacite}
\usepackage{amsfonts}
\usepackage{amssymb}
\usepackage{fancyhdr}
\usepackage{xcolor,graphicx}
\usepackage{setspace}
\usepackage{lscape}
\usepackage{textcomp}
\usepackage[body={6.5in, 9.1in},left=1in,right=1in,top=1in]{geometry}
\usepackage{epsfig}
\usepackage{array}
\usepackage{bbm}
\usepackage{subfigure}
\usepackage{multirow}
\usepackage{amsthm}
\usepackage{latexsym}
\usepackage{epic}
\usepackage{mathrsfs}

\definecolor{light-blue}{rgb}{0.1,0.03,1}
\definecolor{mred}{rgb}{1,0.0,.3}
\begin{document}

\begin{center}
\large{\textbf{Improved Likelihood Estimation for the Generalized Extreme Value and the Inverse Gaussian Lifetime Distributions}}
\end{center}
		\vspace{0.4cm}

\begin{center}
\textbf{Md. Mazharul Islam} \\Applied Statistics\\ Institute of Statistical Research and Training\\ University of Dhaka, Dhaka 1000, Bangladesh
		\vspace{0.2cm}

and\\
		\vspace{0.2cm}

\textbf{ Md Hasinur Rahaman Khan\footnote{corresponding author}} \\ Applied Statistics \\
Institute of Statistical Research and Training \\
University of Dhaka, Dhaka 1000, Bangladesh\\ Email: hasinur@isrt.ac.bd
\end{center}
\begin{center}
\textbf{Abstract}
\end{center}
In presence of nuisance parameters, profile likelihood inference is often unreliable and biased, particularly in small sample scenario.
Over past decades several adjustments have been proposed to modify profile likelihood function in literature including a
modified profile likelihood estimation technique introduced in Barndorff--Nielsen. In this study, adjustment of profile
likelihood function of parameter of interest in presence of nuisance parameter is investigated. We particularly focuss
to extend the Barndorff--Nielsen's technique on Inverse Gaussian distribution for estimating its dispersion parameter and on generalized extreme value (GEV) distribution for estimating its shape parameter. The accelerated failure time models are used for lifetimes having GEV distribution and the Inverse Gaussian distribution is used for lifetime distribution. Monte-Carlo simulation studies are conducted to demonstrate the performances of both approaches. Simulation results suggest the superiority of the modified profile likelihood estimates over the profile likelihood estimates for the parameters of interest. Particularly, it is found that the modifications can improve the overall performance of the estimators through reducing their biases and standard errors.\\
\indent{\textbf{keywords:}} Modified profile likelihood; Profile likelihood; Generalized Extreme Value Distribution; Inverse Gaussian.\\

\newpage
\section{Introduction}
The inverse Gaussian distribution was introduced by Schroedinger (1915). Since this distribution can be regarded as the first passage time in a Brownian motion, it has applications in fields such as economics, biology, medicine and reliability testing [\cite{chhikara1988inverse}, \cite{hung2000information}]. The probability density of an inverse Gaussian distribution is of the form
    \begin{equation}
    f(x;\mu,\lambda) = (\lambda / 2\pi x^{3})^{1/2}  \text{exp}\left\{-\frac{\lambda(x-\mu)^{2}}{2\mu^{2}x}\right\},
    \end{equation}
where $\mu$ is the mean parameter and $\lambda$ is the dispersion parameter, $x >0$, $\mu>0$ and $\lambda>0$. We denote this distribution by $\text{IG}(\mu,\lambda)$, where $\lambda$ is the parameter of interest and $\mu$ is treated as a nuisance parameter [\cite{hung2000information}]. The generalized extreme value (GEV) distribution is considered to be a family of continuous probability distributions that consists of three families--the Gumbel, Frchet, and Weibull which are the extreme value distributions in statistics. The distribution has the density function
     \begin{equation*}
     f(y;\mu,\sigma,\xi) = \frac{1}{\sigma}\left\{1+\xi\left(\frac{y-\mu}{\sigma} \right)  \right\}^{-1/\xi-1}\exp\left[ - \left\{1+\xi \left(\frac{y-\mu}{\sigma}  \right)   \right\}^{-1/\xi} \right],
     \end{equation*}
with the distribution function
    \begin{equation*}
    F(y;\mu,\sigma,\xi) =  \exp\left[ - \left\{1+\xi \left(\frac{y-\mu}{\sigma}  \right)   \right\}^{-1/\xi} \right],
    \end{equation*}
where $\mu\in\Re$ is the location parameter, $\sigma>0$ is the scale parameter and $\xi\in\Re$ is the shape parameter. The shape parameter $\xi$ leads to several distributions, such as Gumbel when $\xi = 0$, Weibull when $\xi < 0$, and Frechet when $\xi > 0$.

The profile likelihood is the most computationally convenient to make inference about a parameter of interest in presence of multidimensional parameter. This is also very widely used procedure of estimating parameters in presence of multidimensional parameter. The technique involves replacing the nuisance parameters in the likelihood function by their maximum likelihood estimates and then examining the resulting profile likelihood as a function of parameter of interest. The profile likelihood brings inconsistent estimates when there is usually a larger number of nuisance parameters [\cite{barndorff1980conditionality}, \cite{cox1987parameter}].
Not many adjustments have been proposed to modify profile likelihood function when there is large number of nuisance parameters. A modified profile likelihood estimation technique introduced in Barndorff--Nielsen (1983)\nocite{barndorff1983formula} is the widely known adjustment among the adjustments. The authors provided a formula which was a synthesis and extension of various results as found in many research studies
including \cite{fisher1934two}, \cite{fraser1968structure}, \cite{daniels1954saddlepoint}, \cite{barndorff1979edgeworth}, \cite{cox1980local},
\cite{hinkley1980likelihood}. The formula leads to a modification of traditional profile likelihood approach when multidimensional parameters exist. In this paper, we particularly have focussed to extend the Barndorff--Nielsen's technique on the dispersion parameter $\lambda$ of inverse Gaussian distribution and also on the shape parameter $\xi$ of generalized extreme value distribution. A similar study has recently been done in \cite{Isla:Khan2015}, where modified profile likelihood is investigated for Weibull regression models for estimating their shape parameters in presence of many nuisance parameters and model regression parameters in presence of collinearity among covariates.

The rest of the manuscript is organized as follows. In Section 2, profile and modified profile likelihood estimation procedures are discussed along with an approximation method to modified profile approach. In Section 3, profile and modified profile likelihood functions corresponding to the dispersion parameter $\lambda$ of inverse Gaussian distribution and shape parameter $\xi$ of Generalized extreme value distribution are derived. In Section 4, several Monte--Carlo simulation results are presented for comparing the overall performances of profile and modified profile likelihood estimators in terms of several statistical measures namely, mean, variance, bias, mean square error and relative bias. The last section discusses the concluding remarks of this study.

\section{Methods}
Suppose, $y_{1}, y_{2},\ldots, y_{n} $ be a sample of $n$ observations from a density $f(y;\,\theta)$, where $\theta$ can be partitioned as $\theta = (\psi,\,\chi)$. Let $\psi$ be the parameter of interest and $\chi$ is considered as a nuisance parameter. Let $\hat{\chi}_{\psi}$ denotes the maximum likelihood estimate of $\chi$ for a given value of $\psi$. The profile likelihood function of $\psi$ is defined as $L_{p}(\psi) = L(\psi,\,\hat{\chi}_{\psi})$. The modified profile likelihood function $L_{mp}(\psi)$ for a parameter of interest $\psi$ with nuisance parameter $\chi$ is defined by
 \begin{equation}
  L_{mp}(\psi) = M(\psi)\, L_{p}(\psi),
  \label{mpfun}
  \end{equation}
where $M$ is a modifying factor and that can be obtained by $M(\psi) = \left|\dfrac{\partial \hat{\chi}}{\partial \hat{\chi}_{\psi} }\right| \mid\hat{j}_{\psi}\mid^{-\frac{1}{2}}$. Here $\mid.\mid$ denotes the absolute value of a matrix determinant, $\dfrac{\partial \hat{\chi}}{\partial \hat{\chi}_{\psi} }$ is the partial derivatives, and $\hat{j}_\psi = j_{\chi\chi}(\psi,\hat{\chi}_{\psi})$ is the observed information on $\chi$ assuming $\psi$ is known \cite{young2005essentials}. In order to obtain modified profile likelihood function one needs to obtain $\left|\dfrac{\partial \hat{\chi}}{\partial \hat{\chi}_{\psi} }\right|$. In many cases this partial derivative does not have any closer form. An alternative expression for $L_{mp}(\psi)$ that does not involve $\dfrac{\partial \hat{\chi}}{\partial \hat{\chi}_{\psi} }$ is available. However, it involves a sample space derivative [\cite{severini2001approximation}, \cite{severini1998approximation}] of the log-likelihood function and the specification of ancillary $a$ such that $(\hat{\psi}, \hat{\chi},a)$ is a minimal sufficient statistic.

An alternative expression of $\dfrac{\partial \hat{\chi}}{\partial \hat{\chi}_{\psi} }$ can be expressed as
\begin{equation}
\dfrac{\partial \hat{\chi}}{\partial \hat{\chi}_{\psi} } =j_{\chi \chi}(\hat{\chi}_{\psi},\,\psi;\, \hat{\chi},\,\hat{\psi},\,a) \ell_{\chi;\,\hat{\chi}}(\hat{\chi}_{\psi},\psi; \hat{\chi},\hat{\psi},a)^{-1},
\label{samd1}
\end{equation}
where
 \begin{equation}
 \ell_{\chi;\,\hat{\chi}}(\hat{\chi}_{\psi},\,\psi;\, \hat{\chi},\,\hat{\psi},\,a) = \frac{\partial}{\partial \hat{\chi}} \left(\frac{\partial \ell(\hat{\chi}_{\psi},\,\psi;\, \hat{\chi},\,\hat{\psi},\,a)}{\partial \chi} \right).
 \label{ell_chi}
 \end{equation}
Here $\ell(\hat{\chi}_{\psi},\,\psi;\, \hat{\chi},\,\hat{\psi},\,a)$ and $j_{\chi \chi}(\hat{\chi}_{\psi},\,\psi;\, \hat{\chi},\,\hat{\psi},\,a)$ are the log likelihood function and the observed information for $\chi$ respectively. They depend on the data only through the minimal sufficient statistic [\cite{da2008improved}]. An alternative formula for \eqref{ell_chi}, which is obtained through approximating ancillary statistic is given by
\begin{equation}
\ell_{\chi;\,\hat{\chi}}(\hat{\chi}_{\psi},\,\psi;\, \hat{\chi},\,\hat{\psi},\,a) = \ell_{\chi;\,y}(\hat{\chi}_{\psi},\,\psi) \hat{V}_\chi,
\label{samd2}
\end{equation}
where $\ell_{\chi,\,y}(\chi,\,\psi) = \frac{\partial \ell_\chi(\chi,\,\psi)}{\partial y}$ and $\hat{V}_\chi = \left(-\frac{\partial F_1(y_1;\,\hat{\chi},\,\hat{\psi})/\partial\hat{\chi}}{p_1(y_1,\,\hat{\chi},\,\hat{\psi})}  \ldots -\frac{\partial F_n(y_n;\,\hat{\chi},\,\hat{\psi})/\partial\hat{\chi}}{p_n(y_n,\hat{\chi},\,\hat{\psi})} \right)^\intercal$ \cite{young2005essentials}. So, the resulting modified profile likelihood function from \eqref{mpfun} takes the form
$$L_{mp}(\psi) = L_p(\psi) |j_{\chi \chi}(\hat{\chi}_{\psi},\,\psi)|^{1/2}|\ell_{\chi;\,\hat{\chi}}(\hat{\chi}_{\psi},\,\psi)|^{-1}$$.
\subsection{Dispersion Parameter of Inverse Gaussian Distribution \label{lkINV}}
Let $x_{1}, x_{2}, \ldots, x_{n}$  be a random sample of size $n$ from $\text{IG}(\mu,\,\lambda)$. The resulting likelihood function and log likelihood function are defined by
   \begin{equation*}
   L(\mu ,\, \lambda) = \left(\frac{\lambda}{2\pi}\right)^{\frac{n}{2}}\left(\prod_{i=1}^{n}x_{i}^{-\frac{3}{2}}\right)
   \text{ exp}\left\{-\frac{\lambda}{2\mu^{2}}\sum_{i=1}^{n}\frac{(x_{i}-\mu)^{2}}{x_{i}}\right\}
   \end{equation*}
and
    \begin{equation}
    \ell(\mu , \lambda) \propto \frac{n}{2}\text{log}\lambda - \frac{\lambda}{2\mu^{2}}
    \sum_{i=1}^{n}\frac{(x_{i}-\mu)^{2}}{x_{i}}
    \label{iglog}
    \end{equation}
respectively. By differentiating equation \eqref{iglog} with respect to $\mu$ and equating to $0$, we get
$\hat{\mu} = \hat{\mu_{\lambda}} = \bar{x}$. The observed information is obtained as
  \begin{align*}
  -\frac{\partial^2}{\partial\mu^2} \ell(\mu,\lambda) & =  \frac{n\lambda}{\bar{x}^3}.
  \end{align*}
The observed information matrix takes the following form $\mid \hat{j_{\lambda}} \mid = \frac{n\lambda}{\bar{x}^{3}}$. The log profile and the log modified profile likelihood functions of $\lambda$ denoted by $\ell_{p}(\lambda)$ and $\ell_{mp}(\lambda)$, respectively, are given by
    \begin{equation*}
    \ell_{p}(\lambda) = \frac{n}{2}\log \lambda -  \frac{\lambda}{2\bar{x}^{2}}
     \sum_{i=1}^{n}\frac{(x_{i}-\bar{x})^{2}}{x_{i}},
    \end{equation*}
    \begin{equation*}
    \ell_{mp}(\lambda) = -\frac{1}{2} \log\frac{n\lambda}{\bar{x}^{3}} + \frac{n}{2}\log\lambda -  \frac{\lambda}{2\bar{x}^{2}}
     \sum_{i=1}^{n}\frac{(x_{i}-\bar{x})^{2}}{x_{i}}.
     \end{equation*}

\subsection{Shape Parameter of Generalized Extreme Value Distribution} \label{lkGEV}
We assume that there exists a number of censored observation in the data. The censoring considered here is assumed to follow a right censoring mechanism. For simulation studies later in this paper, we have assumed that there are 25\% censored observations in the data. The general formula of likelihood function for $n$ random sample from any lifetime distribution having density function $f(t)$ and survivor function $S(t)$ given that some observations are censored, is given by
 \begin{equation}
 L = \prod_{i=1}^{n}\left\{f(t_i) \right\}^{\delta_i} \left\{S(t_i) \right\}^{1-\delta_i},
 \label{lkcen}
 \end{equation}
where $\delta_i = \text{I}(T_i \leq C_i)$ and where $T_i$ and $C_i$ are the $i$-th failure time and censoring time, respectively having independent distribution. Survivor function, $S(t)$ is obtained by
 \begin{align*}
 S(y;\mu,\sigma,\xi) = 1 -   \exp\left[ - \left\{1+\xi \left(\frac{y-\mu}{\sigma}  \right)   \right\}^{-1/\xi} \right].
 \end{align*}

\subsubsection{Profile Likelihood of GEV Regression Model}
Suppose $\mathbf{y} = y_1,\,y_2,\,\ldots,\,y_n$ be $n \times 1$ vector of independent observations so that $y_j \sim GEV(\eta(x_j),\,\sigma,\,\xi)$ for $j=1,\,2,\,\ldots,\,n$. Here $\eta(x_j)= \mathbf{x}_j\phi$, $\mathbf{x}_j = (x_{j1},\, x_{j2},\, \ldots,\, x_{jp})$ are the given set of covariates, and $\phi = (\phi_1,\, \phi_2, \ldots, \phi_p)^\top$. Let $C$ and $\bar{C}$ denote the sets of censored and uncensored observations, respectively. It is noted that in this study, $\xi$ is the parameter of interest, and $\phi$ and $\sigma$ are the nuisance parameters. Assume that $z_j = \frac{y_j-\mathbf{x}_j\phi}{\sigma}$ and $m_j = 1+\xi z_j$. Hence, the density and survivor function will take the form
  \begin{align}
  f(y;\phi,\sigma,\xi) & = \frac{1}{\sigma}m_j^{-1/\xi-1}\exp\left(-m_j^{-1/\xi}\right) , \label{gevfs1} \\
  S(y;\phi,\sigma,\xi) & = 1- \exp\left(-m_j^{-1/\xi}  \right).
  \label{gevfs2}
 \end{align}
Now using \eqref{lkcen}, \eqref{gevfs1} and \eqref{gevfs2} the profile likelihood function will take the form as defined by
  \begin{equation*}
  L_p(\phi,\sigma,\xi) = \prod_{j=1}^{n} \left\{\frac{1}{\sigma}m_j^{-1/\xi-1} \exp\left(-m_j^{-1/\xi}\right) \right\}^{\delta_j} \left\{1-\exp\left(-m_j^{-1/\xi}\right) \right\}^{1-\delta_j}.
  \end{equation*}
The profile log-likelihood, with $r=\sum_{j=1}^{n}\delta_j$, then becomes
    \begin{align*}
   \ell_p(\phi,\sigma,\xi) & = \sum_{j=1}^{n}\delta_j \left\{ \log\frac{1}{\sigma} - (1/\xi+1)\log m_j - m_j^{-1/\xi}\right\} +  (1-\delta_j)\log \left\{1- \exp\left(-m_j^{-1/\xi}\right) \right\}  \\
     & = \sum_{j\in C} \log \left\{1- \exp\left(-m_j^{-1/\xi}\right) \right\} - r \log \sigma - (1/\xi+1)\sum_{j \in \bar{C}} \log m_j - \sum_{j \in \bar{C}} m_j^{-1/\xi}.
     \end{align*}

\subsubsection{Modified Profile Likelihood Function of GEV Distribution}
For obtaining modified version of the profile likelihood function, we have determined the observed information matrix $j_{\phi\sigma}$ and $\ell_{\phi,\,\sigma;\,\hat{\phi},\,\hat{\sigma}}$. The components that are necessary to obtain the modified version of the profile likelihood function can be derived as below.
    \begin{align}
   \frac{\partial m_j}{\partial \phi_s}  = \frac{\partial}{\partial \phi_s} \left\{1+\xi \left(\frac{y_j-x_j\phi}{\sigma}\right)\right\} = \frac{-\xi x_{js}}{\sigma},
   \label{com1}
     \end{align}

  \begin{equation}
   \frac{\partial m_j}{\partial \sigma} =  \frac{\partial}{\partial \sigma} \left\{1+\xi \left(\frac{y_j-x_j\phi}{\sigma}\right)\right\} = \frac{-\xi  z_j}{\sigma},
   \label{com2}
  \end{equation}

  \begin{equation*}
   \frac{\partial r\log\sigma}{\partial m_j}  = r \frac{\partial}{\partial \sigma}\log\sigma \frac{\partial\sigma}{\partial m_j} = -\frac{r}{\xi  z_j}=-\sum_{j=1}^{n}\frac{\delta_j}{\xi z_j}, \\
  \label{com3}
  \end{equation*}

  \begin{equation*}
  \frac{\partial}{\partial m_j}\frac{1}{\sigma} = \frac{\partial}{\partial \sigma}\frac{1}{\sigma} \frac{\partial\sigma}{\partial m_j}=  -\frac{1}{\sigma^2}\left(-\frac{\sigma}{\xi z_j}\right) = \frac{1}{\sigma\xi z_j}, \\
  \end{equation*}

  \begin{equation*}
  \frac{\partial}{\partial m_j}\left(\frac{r}{z_j}\right) = r \frac{\partial}{\partial\sigma} \frac{1}{z_j} \frac{\partial\sigma}{\partial m_j} = -\sum_{j=1}^{n}\frac{\delta_j\sigma}{\xi z_j(y_j-\mathbf{x}_j\phi)}.
  \label{com5}
  \end{equation*}

Now by differentiating the profile log-likelihood, we get
 \begin{equation}
 \frac{\partial \ell(\phi,\sigma,\xi)}{\partial\phi_s} = \frac{\partial \ell(\phi,\sigma,\xi)}{\partial m_j} \times \frac{\partial m_j}{\partial\phi_s}.
 \label{profgevlk}
 \end{equation}

Using \eqref{com1}, we get
 \begin{align*}
 \frac{\partial \ell(\phi,\sigma,\xi)}{\partial m_j} & = -\sum_{j\in C} \frac{(1/\xi)m_j^{-1/\xi-1}}{\left\{\exp\left(m_j^{-1/\xi}\right)-1\right\}} + \sum_{j=1}^{n}\frac{\delta_j}{\xi z_j} - (1/\xi +1)\sum_{j \in \bar{C}}\frac{1}{m_j} + (1/\xi)\sum_{j\in \bar{C}}m_j^{-1/\xi-1},
 \end{align*}
and using \eqref{profgevlk}, we can write
 \begin{align*}
 \ell_{\phi_s} = \frac{\partial \ell(\phi,\sigma,\xi)}{\partial \phi_s} =\frac{x_{js}}{\sigma}\left[\sum_{j\in C} \frac{m_j^{-1/\xi-1}}{\left\{\exp\left(m_j^{-1/\xi}\right)-1\right\}} - \sum_{j=1}^{n}\frac{\delta_j}{z_j} + (1+\xi)\sum_{j\in\bar{C}}\frac{1}{m_j}-
  \sum_{j\in\bar{C}} m_j^{-1/\xi-1}\right]
 \end{align*}
and
    \begin{align}
   \ell_{\phi_s\phi_t} =  \frac{\partial^2\ell(\phi,\sigma,\xi)}{\partial\phi_s\partial\phi_t} & = \frac{\partial\ell_{\phi_s}}{\partial\phi_t}  = \frac{\partial\ell_{\phi_s}}{\partial m_j} \frac{\partial m_j}{\partial \phi_t}.
   \label{der2g}
    \end{align}
Now,
    \begin{align*}
    \frac{\partial\ell_{\phi_s}}{\partial m_j} & = \frac{x_{js}}{\sigma}\frac{\partial}{\partial m_j} \times  \left[\sum_{j\in C} \frac{m_j^{-1/\xi-1}}{\left\{\exp\left(m_j^{-1/\xi}\right)-1\right\}} - \sum_{j=1}^{n}\frac{\delta_j}{z_j} + (1+\xi)\sum_{j\in\bar{C}}\frac{1}{m_j} - \sum_{j\in\bar{C}} m_j^{-1/\xi-1}\right] \\ &
    + \left[\sum_{j\in C} \frac{m_j^{-1/\xi-1}}{\left\{\exp\left(m_j^{-1/\xi}\right)-1\right\}} - \sum_{j=1}^{n}\frac{\delta_j}{z_j} + (1+\xi)\sum_{j\in\bar{C}}\frac{1}{m_j} - \sum_{j\in\bar{C}} m_j^{-1/\xi-1}\right] \times  \frac{\partial}{\partial m_j}\left(\frac{x_{js}}{\sigma}\right),
    \end{align*}
and
    \begin{align*}
    \frac{\partial}{\partial m_j}\left[\sum_{j\in C} \frac{m_j^{-1/\xi-1}}{\left\{\exp\left(m_j^{-1/\xi}\right)-1\right\}} - \sum_{j=1}^{n}\frac{\delta_j}{z_j} + (1+\xi)\sum_{j\in\bar{C}}\frac{1}{m_j} - \sum_{j\in\bar{C}} m_j^{-1/\xi-1}\right] & = \\
    \sum_{j\in C} \left[\frac{(1/\xi)\left(m_j^{-1/\xi-1}\right)^2\exp\left(m_j^{-1/\xi}\right)}{\left\{\exp\left(m_j^{-1/\xi}\right)-1\right\}^2} - \frac{(1/\xi+1)m_j^{-1/\xi-2}}{\left\{\exp\left(m_j^{-1/\xi}\right)-1\right\}} \right] + \sum_{j=1}^{n}\frac{\delta_j\sigma}{\xi z_j (y_j-\mathbf{x}_j\phi)} \\
    -(1+\xi)\sum_{j\in \bar{C}}\frac{1}{m_j^2}+(1/\xi+1)\sum_{j\in \bar{C}}m_j^{-1/\xi-2}
       \end{align*}
that implies
    \begin{align*}
      \frac{\partial\ell_{\phi_s}}{\partial m_j} & =\sum_{j\in C} \left[\frac{(1/\xi)\left(m_j^{-1/\xi-1}\right)^2\exp\left(m_j^{-1/\xi}\right)}{\left\{\exp\left(m_j^{-1/\xi}\right)-1\right\}^2} - \frac{(1/\xi+1)m_j^{-1/\xi-2}}{\left\{\exp\left(m_j^{-1/\xi}\right)-1\right\}} \right] \frac{x_{js}}{\sigma}\\
      & +  \left\{ \sum_{j=1}^{n}\frac{\delta_j\sigma}{\xi z_j (y_j-\mathbf{x}_j\phi)}
         -(1+\xi)\sum_{j\in \bar{C}}\frac{1}{m_j^2}+(1/\xi+1)\sum_{j\in \bar{C}}m_j^{-1/\xi-2} \right\} \frac{x_{js}}{\sigma} \\
      & +  \left[\sum_{j\in C} \frac{m_j^{-1/\xi-1}}{\left\{\exp\left(m_j^{-1/\xi}\right)-1\right\}} - \sum_{j=1}^{n}\frac{\delta_j}{z_j} + \sum_{j\in\bar{C}}\frac{1+\xi}{m_j} - \sum_{j\in\bar{C}} m_j^{-1/\xi-1}  \right]\frac{x_{js}}{\sigma\xi z_j}.
    \end{align*}

Using equation \eqref{der2g}, we can write
 \begin{align*}
 \frac{\partial^2\ell(\phi,\sigma,\xi)}{\partial\phi_s\partial\phi_t} & =  \sum_{j\in C} \left[\frac{(1/\xi)\left(m_j^{-1/\xi-1}\right)^2\exp\left(m_j^{-1/\xi}\right)}{\left\{\exp\left(m_j^{-1/\xi}\right)-1\right\}^2} - \frac{(1/\xi+1)m_j^{-1/\xi-2}}{\left\{\exp\left(m_j^{-1/\xi}\right)-1\right\}} \right] \times \frac{x_{js}}{\sigma}\left(-\frac{\xi x_{jt}}{\sigma}\right)\\
 & +  \left\{ \sum_{j=1}^{n}\frac{\delta_j\sigma}{\xi z_j (y_j-\mathbf{x}_j\phi)}
    -(1+\xi)\sum_{j\in \bar{C}}\frac{1}{m_j^2}+(1/\xi+1)\sum_{j\in \bar{C}}m_j^{-1/\xi-2} \right\} \times \frac{x_{js}}{\sigma}\left(-\frac{\xi x_{jt}}{\sigma}\right) \\
 & +  \left[\sum_{j\in C} \frac{m_j^{-1/\xi-1}}{\left\{\exp\left(m_j^{-1/\xi}\right)-1\right\}} - \sum_{j=1}^{n}\frac{\delta_j}{z_j} + \sum_{j\in\bar{C}}\frac{1+\xi}{m_j} - \sum_{j\in\bar{C}} m_j^{-1/\xi-1}  \right] \times \frac{x_{js}}{\sigma\xi z_j}\left(-\frac{\xi x_{jt}}{\sigma}\right).
 \end{align*}

Again using results from \eqref{com2}, we have
 \begin{align*}
 \ell_\sigma & = \frac{\partial\ell(\phi,\sigma,\xi)}{\partial\sigma} = \frac{\partial\ell(\phi,\sigma,\xi)}{\partial m_j} \frac{\partial m_j}{\partial\sigma}\\
 & = \left[\sum_{j\in C} \frac{m_j^{-1/\xi-1}}{\left\{\exp\left(m_j^{-1/\xi}\right)-1\right\}} - \sum_{j=1}^{n}\frac{\delta_j}{z_j} + (1+\xi)\sum_{j\in\bar{C}}\frac{1}{m_j} - \sum_{j\in\bar{C}} m_j^{-1/\xi-1} \right] \frac{z_j}{\sigma},
 \end{align*}
and $\frac{\partial^2\ell(\phi,\sigma,\xi)}{\partial\sigma^2} = \frac{\partial\ell_\sigma}{\partial\sigma} = \frac{\partial}{\partial m_j} \ell_\sigma\frac{\partial m_j}{\partial\sigma}$. Another component is derived as
  \begin{align}
 \frac{\partial}{\partial m_j}\left(\frac{z_j}{\sigma}\right) & = \frac{\partial}{\partial\sigma}\left(\frac{z_j}{\sigma}\right)\frac{\partial\sigma}{\partial m_j} = \frac{\partial}{\partial\sigma}\left(\frac{y_j-\mathbf{x}_j\phi}{\sigma^2}\right)\frac{\partial\sigma}{\partial m_j}= 2 \left(\frac{y_j-\mathbf{x}_j\phi}{\sigma^2}\right) \frac{\sigma}{\xi (y_j-\mathbf{x}_j\phi)} = \frac{2}{\xi \sigma}.
 \label{acom}
  \end{align}

Using \eqref{acom}, we can write
    \begin{align*}
    \frac{\partial\ell_\sigma}{\partial m_j} & = \frac{z_j}{\sigma}\frac{\partial}{\partial m_j} \left[\sum_{j\in C} \frac{m_j^{-1/\xi-1}}{\left\{\exp\left(m_j^{-1/\xi}\right)-1\right\}} - \sum_{j=1}^{n}\frac{\delta_j}{z_j} + (1+\xi)\sum_{j\in\bar{C}}\frac{1}{m_j} - \sum_{j\in\bar{C}} m_j^{-1/\xi-1} \right] \\
    & + \left[\sum_{j\in C} \frac{m_j^{-1/\xi-1}}{\left\{\exp\left(m_j^{-1/\xi}\right)-1\right\}} - \sum_{j=1}^{n}\frac{\delta_j}{z_j} + (1+\xi)\sum_{j\in\bar{C}}\frac{1}{m_j} - \sum_{j\in\bar{C}} m_j^{-1/\xi-1} \right] \frac{\partial}{\partial m_j}\left(\frac{z_j}{\sigma}\right) \\
    & = \left[\sum_{j\in C}\left\{\frac{(1/\xi)\left(m_j^{-1/\xi-1}\right)^2\exp\left(m_j^{-1/\xi}\right)}{\left\{\exp\left(m_j^{-1/\xi} \right)-1 \right\}^2} \right\} - \frac{(1/\xi+1)m_j^{-1/\xi-2} }{\left\{\exp\left(m_j^{-1/\xi} \right)-1 \right\}} \right] \frac{z_j}{\sigma}\\
    & +  \left\{\sum_{j=1}^{n}\frac{\delta_j\sigma}{\xi z_j (y_j-\mathbf{x}_j\phi)} - (1+\xi)\sum_{j \in \bar{C}}\frac{1}{m_j^2}+(1/\xi+1)\sum_{j\in \bar{C}}m_j^{-1/\xi-1}  \right\} \frac{z_j}{\sigma} \\
    & + \left[ \sum_{j\in C} \frac{m_j^{-1/\xi-1}}{\left\{\exp\left(m_j^{-1/\xi}\right)-1\right\}} - \sum_{j=1}^{n}\frac{\delta_j}{z_j} + (1+\xi)\sum_{j\in\bar{C}}\frac{1}{m_j} - \sum_{j\in\bar{C}} m_j^{-1/\xi-1} \right] \frac{2}{\xi\sigma}
    \end{align*}
and so,
    \begin{align*}
    \frac{\partial^2\ell(\phi,\sigma,\xi)}{\partial\sigma^2} & = \left[\sum_{j\in C}\left\{\frac{(1/\xi)\left(m_j^{-1/\xi-1}\right)^2\exp\left(m_j^{-1/\xi}\right)}{\left\{\exp\left(m_j^{-1/\xi} \right)-1 \right\}^2} \right\} - \frac{(1/\xi+1)m_j^{-1/\xi-2} }{\left\{\exp\left(m_j^{-1/\xi} \right)-1 \right\}} \right] \times \left(-\frac{\xi z_j}{\sigma}\right)  \frac{z_j}{\sigma}\\
    & +  \left\{\sum_{j=1}^{n}\frac{\delta_j\sigma}{\xi z_j (y_j-\mathbf{x}_j\phi)} - (1+\xi)\sum_{j \in \bar{C}}\frac{1}{m_j^2}+(1/\xi+1)\sum_{j\in \bar{C}}m_j^{-1/\xi-1}  \right\} \times \left(-\frac{\xi z_j}{\sigma}\right) \frac{z_j}{\sigma} \\
    & + \left[ \sum_{j\in C} \frac{m_j^{-1/\xi-1}}{\left\{\exp\left(m_j^{-1/\xi}\right)-1\right\}} - \sum_{j=1}^{n}\frac{\delta_j}{z_j} + (1+\xi)\sum_{j\in\bar{C}}\frac{1}{m_j} - \sum_{j\in\bar{C}} m_j^{-1/\xi-1} \right] \times \left(-\frac{\xi z_j}{\sigma}\right) \frac{2}{\xi\sigma}.
       \end{align*}
Now,
 \begin{align*}
 \frac{\partial^2\ell(\phi,\sigma\xi)}{\partial\sigma\partial\phi_s} = \frac{\partial}{\partial\sigma}\ell_{\phi_s} = \frac{\partial}{\partial m_j} \ell_{\phi_s} \frac{\partial m_j}{\partial\sigma}
 \end{align*}
which takes the following form by using result of \eqref{com2}.
 \begin{align*}
 \frac{\partial^2\ell(\phi,\sigma\xi)}{\partial\sigma\partial\phi_s} & = \sum_{j\in C} \left[\frac{(1/\xi)\left(m_j^{-1/\xi-1}\right)^2\exp\left(m_j^{-1/\xi}\right)}{\left\{\exp\left(m_j^{-1/\xi}\right)-1\right\}^2} - \frac{(1/\xi+1)m_j^{-1/\xi-2}}{\left\{\exp\left(m_j^{-1/\xi}\right)-1\right\}} \right] \times \frac{x_{js}}{\sigma}\left(-\frac{\xi z_{j}}{\sigma}\right)\\
 & +  \left\{ \sum_{j=1}^{n}\frac{\delta_j\sigma}{\xi z_j (y_j-\mathbf{x}_j\phi)}
    -(1+\xi)\sum_{j\in \bar{C}}\frac{1}{m_j^2}+(1/\xi+1)\sum_{j\in \bar{C}}m_j^{-1/\xi-2} \right\} \times  \frac{x_{js}}{\sigma}\left(-\frac{\xi z_{j}}{\sigma}\right) \\
 & +  \left[\sum_{j\in C} \frac{m_j^{-1/\xi-1}}{\left\{\exp\left(m_j^{-1/\xi}\right)-1\right\}} - \sum_{j=1}^{n}\frac{\delta_j}{z_j} + \sum_{j\in\bar{C}}\frac{1+\xi}{m_j} - \sum_{j\in\bar{C}} m_j^{-1/\xi-1}  \right] \times  \frac{x_{js}}{\sigma\xi z_j}\left(-\frac{\xi z_{j}}{\sigma}\right).
 \end{align*}

The $(p+1)\times(p+1)$ observed information matrix now can be obtained using all the components as derived before. Now, we need to derive an expression for $\left|\frac{\partial\hat{\chi}}{\partial\hat{\chi}_\xi} \right|$, where $\chi=(\phi_1,\,\phi_2,\,\ldots,\,\phi_p,\,\sigma)$. Using equations \eqref{samd1}, \eqref{samd2} and \eqref{ell_chi}, an alternative approximation can be obtained via sample space derivatives. This calculation involves
\begin{align*}
 \ell_{\chi;\hat{\chi}}(\xi) = \frac{\partial}{\partial\hat{\chi}} \ell_\chi(\xi) = \frac{\partial}{\partial y}\ell_\chi(\xi) \frac{\partial y}{\partial\hat{\chi}} = \frac{\partial}{\partial m_j}\ell_\chi(\xi) \frac{\partial m_j}{\partial y} \mathbf{\hat{V}_\chi},
 \end{align*}
where $\frac{\partial m_j}{\partial y} = \frac{\partial}{\partial y} \left\{ 1+\xi\left(\frac{y_j-\mathbf{x}_j\phi}{\sigma} \right) \right\} = \frac{\xi}{\sigma}$ and $\hat{V}_\chi = \left(-\frac{\partial F_1(y_1;\hat{\chi},\hat{\psi})/\partial\hat{\chi}}{p_1(y_1,\hat{\chi},\hat{\psi})}  \ldots -\frac{\partial F_n(y_n;\hat{\chi},\hat{\psi})/\partial\hat{\chi}}{p_n(y_n,\hat{\chi},\hat{\psi})} \right)^\intercal $. The following component has already been derived earlier.
  \begin{equation*}
  \frac{\partial}{\partial m_j}\ell_\chi(\xi) = \left\{\frac{\partial}{\partial m_j}\ell_{\phi_1}(\xi), \frac{\partial}{\partial m_j}\ell_{\phi_2}(\xi),\ldots, \frac{\partial}{\partial m_j}\ell_{\phi_p}(\xi), \frac{\partial}{\partial m_j}\ell_{\sigma}(\xi) \right\}.
  \end{equation*}
Now,
 \begin{align*}
 \frac{\partial}{\partial\phi}F(y;\phi,\sigma,\xi) & = \frac{\partial}{\partial\phi}\exp\left\{-m_j^{-1/\xi} \right\} = \frac{\partial}{\partial m_j}\exp\left\{-m_j^{-1/\xi} \right\} \frac{\partial m_j}{\partial\phi} \\
 & = \exp\left\{-m_j^{-1/\xi} \right\} \left(\frac{1}{\xi} m_j^{-1/\xi-1} \right) \left( - \frac{\xi x_j}{\sigma} \right) \\
 & = -m_j^{-1/\xi-1} \exp\left(-m_j^{-1/\xi} \right) \frac{x_j}{\sigma}
 \end{align*}
and hence
 \begin{align*}
 \frac{\frac{\partial}{\partial\phi}F(y;\phi,\sigma,\xi)}{f(y;\phi,\sigma,\xi)} = \frac{-m_j^{-1/\xi-1}\exp\left( -m_j^{-1/\xi} \right)\frac{x_j}{\sigma}}{\frac{1}{\sigma}m_j^{-1/\xi-1}\exp\left( -m_j^{-1/\xi}\right)} = -x_j.
 \end{align*}
 This, we have
 \begin{align*}
 \frac{\partial}{\partial\sigma}F(y;\phi,\sigma,\xi) & = \frac{\partial}{\partial\phi}\exp\left\{-m_j^{-1/\xi} \right\} = \frac{\partial}{\partial m_j}\exp\left\{-m_j^{-1/\xi} \right\} \frac{\partial m_j}{\partial\sigma} \\
 & = \exp\left\{-m_j^{-1/\xi} \right\} \left(\frac{1}{\xi} m_j^{-1/\xi-1} \right) \left( - \frac{\xi z_j}{\sigma} \right) \\
 & = -m_j^{-1/\xi-1} \exp\left(-m_j^{-1/\xi} \right) \frac{z_j}{\sigma}
 \end{align*}
and hence
  \begin{align*}
  \frac{\frac{\partial}{\partial\sigma}F(y;\phi,\sigma,\xi)}{f(y;\phi,\sigma,\xi)} = \frac{-m_j^{-1/\xi-1}\exp\left( -m_j^{-1/\xi} \right)\frac{z_j}{\sigma}}{\frac{1}{\sigma}m_j^{-1/\xi-1}\exp\left( -m_j^{-1/\xi}\right)} = -z_j.
  \end{align*}

The above function can be attained easily for uncensored observations but for censored observations the function will take values $\mathbf{0}=(0,\,0,\,\ldots,0)$. Now, we can easily find $\ell_{\chi;\,\hat{\chi}}(\xi)$. Hence, by accumulating all these components, any one can obtain the modified version of profile likelihood function for the shape parameter $\xi$ of generalized extreme value distribution.

\section{Numerical Analysis}
We have conducted several Monte-Carlo simulation studies depending on various sizes of sample considered in the study. The lifetime and censoring times are generated from two independent distributions in a way so that certain censoring rate is maintained. The random right censoring mechanism is adopted in generating observations to be used for the generalized extreme value regression model. The simulations are based on 1000 runs and different summary statistics namely, mean, variance, bias, mean squared error (MSE), and relative bias (RB) are used for comparison purpose. The relative bias is calculated using the formula \texttt{bias/true parameter} and expressed as percentage (\cite{da2008improved}).

\subsection{Example I}
The parameter of interest, $\lambda$ is set to $4$ and the nuisance parameter, $\mu$ is set to $2$. The sample sizes that range from 3 to 50 are considered for the example. Of interest, very low, low medium and low moderated samples are considered. The lifetimes considered are assumed to be all uncensored.

\begin{table}[t]
   \centering
   \caption{Simulation Results of Inverse Gaussian Dispersion Parameter $\lambda = 4$.}
   \scalebox{.8}{\begin{tabular}{|c|cc|cc|cc|cc|cc|}
     \hline
    Sample & \multicolumn{2}{c|}{mean} & \multicolumn{2}{c|}{variance} & \multicolumn{2}{c|}{bias} & \multicolumn{2}{c|}{MSE} &\multicolumn{2}{c|}{RB (\%)}  \\ \cline{2-11}
     size& p & mp & p & mp & p & mp & p & mp & p & mp \\
     \hline
    3 & 8.799 & 7.316 & 22.134 & 22.581 & 4.799 & 3.316 & 45.165 & 33.580 & 119.976 & 82.910 \\
      5 & 7.181 & 6.164 & 16.904 & 15.640 & 3.181 & 2.164 & 27.022 & 20.324 & 79.522 & 54.106 \\
      7 & 6.363 & 5.624 & 12.854 & 11.673 & 2.363 & 1.624 & 18.440 & 14.309 & 59.086 & 40.596 \\
      9 & 5.609 & 5.044 & 8.881 & 7.749 & 1.609 & 1.044 & 11.469 & 8.838 & 40.223 & 26.088 \\
      11 & 5.286 & 4.835 & 6.550 & 5.835 & 1.286 & 0.835 & 8.204 & 6.533 & 32.153 & 20.882 \\
      13 & 5.121 & 4.742 & 5.477 & 4.847 & 1.121 & 0.742 & 6.735 & 5.398 & 28.033 & 18.551 \\
      15 & 4.966 & 4.644 & 4.168 & 3.760 & 0.966 & 0.644 & 5.101 & 4.175 & 24.150 & 16.107 \\
      17 & 4.843 & 4.559 & 3.219 & 2.845 & 0.843 & 0.559 & 3.929 & 3.158 & 21.078 & 13.984 \\
      19 & 4.730 & 4.481 & 3.133 & 2.812 & 0.730 & 0.481 & 3.666 & 3.044 & 18.254 & 12.037 \\
      25 & 4.471 & 4.293 & 1.740 & 1.602 & 0.471 & 0.293 & 1.962 & 1.688 & 11.783 & 7.316 \\
      30 & 4.377 & 4.231 & 1.513 & 1.413 & 0.377 & 0.231 & 1.655 & 1.466 & 9.418 & 5.774 \\
      50 & 4.313 & 4.227 & 0.789 & 0.758 & 0.313 & 0.227 & 0.887 & 0.809 & 7.822 & 5.665 \\
      \hline
   \end{tabular}\label{tab:1}
   }\\ \footnotesize{Here p = profile and mp = modified profile}
   \end{table}

Table \ref{tab:1} presents the Monte-Carlo simulation results of both type of estimators (modified profile and only profile) for the inverse Gaussian dispersion parameter $\lambda$ as discussed in Section \ref{lkINV}. We notice that the modified estimators provide more accurate average maximum likelihood estimates. Typically, the increasing of sample size reduces the MSE and variances. This relationship is expected, because both estimators are asymptotically unbiased. As the variances under profile likelihood are much higher than that of the modified version, the wider confidence intervals for the estimates under profile likelihood estimation to be obtained. It reveals from the results that the modified estimators are the best performer in terms of all statistics under consideration. Particularly, when sample size is low the modified estimators significantly outperform the profile estimators. This result provides even stronger evidence in favor of the accuracy of the modified estimators. The modified estimators consistently exhibit the smallest (35-40\% less) relative bias. Hence, modified estimators generally perform well enduring the variations of the sample sizes.

\subsection{Example II}
We are interested to focus on how modified profile likelihood estimates of $\xi$ differ from its profile likelihood estimates when there is a set of covariates in the model under generalized extreme value AFT models. The profile and modified profile likelihood estimation techniques, under this situation, are discussed in Section \ref{lkGEV}. Of interest, we assume that the true value of the shape parameter is $\xi=2$. A set of two covariates $(x_1,\,x_2)$ is considered and we assume $x_1=1$ and $x_2$ is assumed to be generated as independent realizations from Uniform distribution. The parameters are set as $\sigma=1$, $\phi_1=1$ and $\phi_2=1$ and let $\eta(\mathbf{x}) = \phi_1\,x_1+\phi_2\,x_2$ in consistent with the notations as discussed in Section \ref{lkGEV}. The lifetimes and censoring times are generated from two distributions in a way so that 25\% censoring rate is maintained. A low and a moderated size of sample ($n=20,\, 50$) are considered for this example. The results are presented in Table \ref{Shape parameter  of GEV regression}.

    \begin{table}[h]
  \caption{Simulation Results of Shape Parameter $\xi=2$ of GEV Regression.}
    \centering
  \scalebox{.8}{ \begin{tabular}{|c|cc|cc|cc|cc|cc|}
      \hline
      Sample & \multicolumn{2}{c|}{mean} & \multicolumn{2}{c|}{variance} & \multicolumn{2}{c|}{bias} & \multicolumn{2}{c|}{MSE} & \multicolumn{2}{c|}{RB (\%)} \\ \cline{2-11}
      size & p & mp & p & mp & p & mp & p & mp & p & mp \\
      \hline
     20 & 2.7693 & 2.4600 & 1.8669 & 0.7288 & 0.7693 & 0.4600 & 2.4587 & 0.9405 & 38.4651 & 23.0018 \\
       50 & 2.7507 & 2.7114 & 1.6646 & 0.1114 & 0.7507 & 0.7114 & 2.2281 & 0.6174 & 37.5348 & 35.5690 \\
       \hline
    \end{tabular}}
    \label{Shape parameter  of GEV regression}\\
      \footnotesize{Here p = profile and mp = modified profile}
    \end{table}

As the previous example, the modified estimators consistently exhibit the consistent maximum likelihood average estimates when the sample size considered is very small. The modified estimators provide much less variance, bias and MSE as compared with the original profile estimators. However, when sample size is large i.e. $n=50$ the differences of the statistics between the approaches are found to be very marginal. This indicates that both estimators become asymptotically unbiased for large sample as expected theoretically.

\section{Discussion and Conclusion}
The maximum likelihood estimator can be considerably biased when inference is made based on small sample. In this paper, we have shown how to extend the Barndorff--Nielsen's (1980, 1983) modified profile likelihood estimation techniques on Inverse Gaussian distribution for estimating its dispersion parameter and on generalized extreme value distribution for estimating its shape parameter. The implementations of both approaches are illustrated in detail with two simulated examples by considering variety of sample sizes that range from very small to medium and the nature of lifetimes--uncensored or censored. It reveals from the results of both examples that the modified profile likelihood estimators outperform the usual profile estimators in terms of all statistics used under all variety of sample sizes considered. Particularly, when sample size is low the modified estimators significantly outperform the profile estimators. The results demonstrate even very stronger evidence in favor of the accuracy of the modified estimators. The modified estimators consistently exhibit the smallest relative bias even for the reasonably large samples. Both estimators tend to be unbiased as sample size increases that supports indirectly the claim that the estimators are asymptotically unbiased.

The parametric inference is some how affected by the achievement of parametric orthogonality. In presence of high dimensional nuisance parameters which is unlike the cases, profile likelihood inference technique may quiet suspicious and sometimes unreliable. But, an advantage of the Barndorff--Nielsen's approach of modification is that it does not require both the interest and nuisance parameters to be orthogonal. Furthermore, the modification doesn't require generally the specification of an ancillary statistic. It is interesting to note that the modification to profile likelihood comes from several higher order approximations. The model assumes that the derivatives of the likelihood components and the information matrix are computationally convenient or at least numerically obtainable. If the model does not belong to the family of exponential distribution, the exact or even approximate expression for the conditional distribution of the maximum likelihood estimator will be accurate to order $O(n^{-1})$ and often up to order $O(n{-3/2})$.
 \bibliography{Thesis2}
\bibliographystyle{apacite}

\end{document}